\documentclass[12pt]{article}
\usepackage{latexsym, amssymb, amsmath, amscd, amsfonts, epsfig, graphicx, colordvi,verbatim,ifpdf}
\usepackage{amsfonts, amsmath, amssymb}
\usepackage{amssymb,amsfonts,amsmath,latexsym,epsfig,cite, psfrag,eepic,color}
\usepackage{amscd,graphics}
\usepackage{latexsym, amssymb,  amsmath,amscd, amsfonts, epsfig, graphicx, colordvi,amsthm}

\usepackage{graphicx}
\usepackage{color}
\usepackage{ifpdf}
\usepackage{fancybox}

\usepackage{float}

\newtheorem{theo}{Theorem}

\def\qed{\nopagebreak\hfill{\rule{4pt}{7pt}}
\medbreak}

\def\qed{\nopagebreak\hfill{\rule{4pt}{7pt}}
\medbreak}

\newlength{\boxedparwidth}
  {\begin{center} \begin{tabular}{|@{\hspace{.315in}}c@{\hspace{.15in}}|}
                  \hline \\ \begin{minipage}[t]{\boxedparwidth}
                  \setlength{\parindent}{.25in}}%
  {\end{minipage} \\ \\ \hline \end{tabular} \end{center}}

\begin{document}

\begin{center}

 {\Large \bf  A Combinatorial Proof of a Schmidt  Type Theorem\\[5pt] of Andrews and Paule}

\end{center}

\begin{center}
{Kathy Q. Ji}  \vskip 2mm

 Center for Applied Mathematics\\[3pt]
   Tianjin University\\[3pt]
    Tianjin 300072, P.R. China\\[6pt]
   \vskip 2mm

    kathyji@tju.edu.cn
\end{center}

\vskip 6mm \noindent {\bf Abstract.} This note is devoted to  a combinatorial proof of a Schmidt type theorem due to Andrews and Paule. A four-variable refinement of Andrews and Paule's theorem is also obtained based on this combinatorial construction.

 \vskip 6mm

The main objective of this note is  to give a combinatorial proof of the following  partition theorem due to Andrews and Paule \cite{Andrews-Paule-2021}. Sylvester's bijection \cite{bes94, bre99,  mac84} for Euler's partition theorem and Wright's bijection \cite{and84,yee03,Wright-1965} for the Jacobi's triple product identity plays an important role in the combinatorial construction.

\begin{theo}[Andrews-Paule]\label{main-AP} Assume that $n\geq 1$.
Let $s(n)$ denote the number of partitions $a_1+a_2+a_3+\cdots$ satisfying $a_1\geq a_2\geq
a_3 \geq \cdots$ and $n=a_1+a_3+a_5+\cdots$. Let $t(n)$ denote the number of two-color partitions of $n$.
Then
\begin{equation*}
 s(n)=t(n)
 \end{equation*}
\end{theo}

For example, let $n=3$. There are ten partitions counted by $s(3)$, which are
\[
\begin{array}{ll}
3, 3 + 3, 3 + 2, 3 + 1, 2 + 2 + 1, 2 + 2 + 1 + 1,
2 + 1 + 1,2 + 1 + 1 + 1,&\\[5pt]
  1 + 1 + 1 + 1 + 1, 1 + 1 + 1 + 1 + 1 + 1. &
\end{array}
\]
and there are also ten red and green partitions counted by $t(3)$, which are
\[\begin{array}{ll}
3_r, 3_g, 2_r+1_r, 2_g+1_r, 2_r+1_g, 2_g+1_g, 1_r+1_r+1_r,
 1_r+1_r+1_g, &\\[5pt] 1_r+1_g+1_g,
1_g+1_g+1_g.
\end{array}\]

\noindent{\it Proof.} Let $\mathcal{T}(n)$ denote the set of two-color partitions counted by $t(n)$ and let $\mathcal{S}(n)$ denote the set of partitions counted by $s(n)$. We aim to construct a bijection $\phi$ between $\mathcal{T}(n)$ and $\mathcal{S}(n)$.

Let $\lambda$ be a two-color partition in $\mathcal{T}(n)$ with $r$ red parts and $l$ green parts. Assume that $m=\max\{r,l\}$. We aim to define $\phi(\lambda)=\gamma=(\gamma_1,\gamma_2,\ldots,\gamma_{2m-1},\gamma_{2m})$ such that $\gamma_1\geq \gamma_2\geq \ldots\geq \gamma_{2m-1}\geq \gamma_{2m}\geq 0$  and $\gamma_1+\gamma_3+\cdots+\gamma_{2m-1}=n$.

Let $\alpha$ be a partition consisting of all red parts in $\lambda$      and  $\beta$ be a partition consisting of all green parts in $\lambda$. First, add $0$ at the end of $\alpha$ or $\beta$ so that they are of the same length depending on which is of smaller length.  Assume that $r\leq l$, so $m=l$. Then  $\alpha=(\alpha_1,\alpha_2,\ldots,\alpha_r,\underbrace{0,\ldots,0}_{l-r})$    and  $\beta=(\beta_1,\beta_2,\ldots,\beta_l)$.

We next define a pair $(\overline{\alpha},\overline{\beta})$ of partitions with distinct parts corresponding to $(\alpha,\beta)$, where
$\overline{\alpha}=(\alpha_1+l-1,\alpha_2+l-2,\ldots,
\alpha_r+l-r, l-r-1,\ldots,1,0),$ and $\overline{\beta}=(\beta_1+l-1,\beta_2+l-2,\ldots,
\beta_{l}).$ Obviously, $|\overline{\alpha}|+|\overline{\beta}|=|\alpha|
+|\beta|+l(l-1)$.

We now apply Wright's bijection to represent $(\overline{\alpha},\overline{\beta})$ as a Young diagram of an ordinary partition $Y(\overline{\alpha},\overline{\beta})$: put $l$ squares on the diagonal, and then  for
$j=1, 2, \ldots, l$, put $\bar{\alpha}_j$ squares in row $j$ to the
right of the diagonal and $\bar{\beta}_j$ squares in column $j$ below
the diagonal. For example,  Figure 1. gives the Young diagram of $(\overline{\alpha},\overline{\beta})$, where $\overline{\alpha}=(3,2,0)$ and $\overline{\beta}=(5,3,1)$.

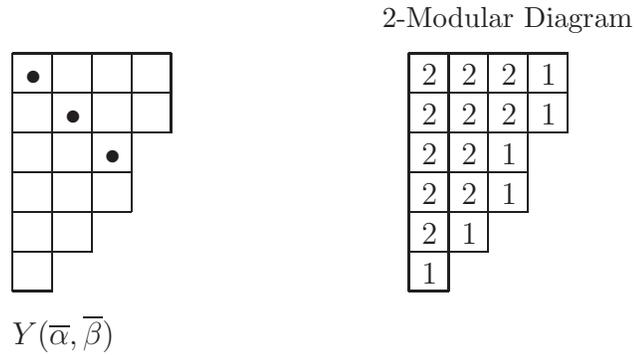
\begin{figure}[h]
\begin{center}
\begin{picture}(100,120)
\put(50,0){\put(-100,100){\line(1,0){60}} \put(-100,100){\line(0,-1){90}}
\put(-100,10){\line(1,0){15}} \put(-85,10){\line(0,1){90}}
\put(-100,25){\line(1,0){30}}\put(-70,25){\line(0,1){75}}
\put(-100,40){\line(1,0){45}}\put(-55,40){\line(0,1){60}}
\put(-100,55){\line(1,0){45}}\put(-100,70){\line(1,0){60}}
\put(-40,70){\line(0,1){30}}\put(-100,85){\line(1,0){60}}\put(-95,88){$\bullet$}
\put(-80,73){$\bullet$} \put(-65,58){$\bullet$}\put(-100,-10){$Y(\overline{\alpha},\overline{\beta})$}
}

\put(200,0){
\put(-100,100){\line(1,0){60}} \put(-100,100){\line(0,-1){90}}
\put(-100,10){\line(1,0){15}} \put(-85,10){\line(0,1){90}}
\put(-100,25){\line(1,0){30}}\put(-70,25){\line(0,1){75}}
\put(-100,40){\line(1,0){45}}\put(-55,40){\line(0,1){60}}
\put(-100,55){\line(1,0){45}}\put(-100,70){\line(1,0){60}}
\put(-40,70){\line(0,1){30}}\put(-100,85){\line(1,0){60}}
\put(-110,110){{\small 2-Modular Diagram}}
\put(-95,88){2}\put(-80,88){2}\put(-65,88){2}\put(-50,88){1}
\put(-95,73){2}\put(-80,73){2}\put(-65,73){2}\put(-50,73){1}
\put(-95,58){2}\put(-80,58){2}\put(-65,58){1} \put(-95,43){2}
\put(-80,43){2}\put(-65,43){1}\put(-95,28){2}\put(-80,28){1}
\put(-95,13){1}}

\end{picture}
\end{center}
\caption{The Young diagram of $(\bar{\alpha},\bar{\beta})$ and 2-Modular diagram.}
\end{figure}

For each row in  the Young diagram of $(\overline{\alpha},\overline{\beta})$, write 2 in each box and a 1 at the end of the row  to obtain the 2-modular diagram. Decompose the 2-modular diagram
into hooks $H_1,\,H_2,\ldots$ with the diagonal boxes as corners.
Let $\mu_1$ be the number of squares in $H_1$, let $\mu_2$ be the
number of 2's in $H_1$, let $\mu_3$ be the number of squares in
$H_2$, let $\mu_4$ be the number of 2's in $H_2$, and so on. Set
$\mu=(\mu_1,\, \mu_2,\,  \ldots,\,\mu_{2l-1},\mu_{2l})$, see
Figure $2$. Then $\mu$ is clearly a partition with distinct parts. Furthermore, $\mu_1+\mu_3+\cdots+\mu_{2l-1}=|\alpha|+|\beta|+l^2$. Hence we may define  $\gamma=(\mu_1-(2l-1),\, \mu_2-{2l-2},\,  \ldots,\,\mu_{2l-1}-1,\mu_{2l}).$ Clearly, $\gamma_1+\gamma_3+\cdots+\gamma_{2l-1}=|\alpha|+|\beta|$, and so $\gamma \in \mathcal{S}(n)$. Furthermore, this process is reversible since Sylvester's bijection and Wright's bijection are reversible. Thus, we complete the proof of Theorem \ref{main-AP}.
\qed

\begin{figure}[h]
\begin{center}
\begin{picture}(100,120)
\put(-100,100){\line(1,0){60}} \put(-100,100){\line(0,-1){90}}
\put(-100,10){\line(1,0){15}} \put(-85,10){\line(0,1){90}}
\put(-100,25){\line(1,0){30}}\put(-70,25){\line(0,1){75}}
\put(-100,40){\line(1,0){45}}\put(-55,40){\line(0,1){60}}
\put(-100,55){\line(1,0){45}}\put(-100,70){\line(1,0){60}}
\put(-40,70){\line(0,1){30}}\put(-100,85){\line(1,0){60}}
\put(-110,110){{\small 2-Modular Diagram}}
\put(-95,88){2}\put(-80,88){2}\put(-65,88){2}\put(-50,88){1}
\put(-95,73){2}\put(-80,73){2}\put(-65,73){2}\put(-50,73){1}
\put(-95,58){2}\put(-80,58){2}\put(-65,58){1} \put(-95,43){2}
\put(-80,43){2}\put(-65,43){1}\put(-95,28){2}\put(-80,28){1}
\put(-95,13){1} {\thicklines
\put(-95,93){\line(1,0){50}}\put(-95,93){\line(0,-1){75}}
\put(-98,15){{\tiny 9}}\put(-90,90){\line(1,0){30}}
\put(-90,90){\line(0,-1){62}}\put(-90,23){{\tiny 7}}
\put(-80,78){\line(1,0){35}}\put(-80,78){\line(0,-1){48}}
\put(-80,25){{\tiny6}}\put(-75,75){\line(1,0){18}}
\put(-75,75){\line(0,-1){32}}\put(-75,39){{\tiny 4}}
\put(-60,65){\line(0,-1){20}}\put(-60,40){{\tiny 2}}
\put(10,50){\vector(1,0){30}}\put(30,50){\vector(-1,0){30}}}
\put(100,100){\line(1,0){135}}\put(100,100){\line(0,-1){75}}
\put(100,25){\line(1,0){30}}\put(130,25){\line(0,1){75}}
\put(100,40){\line(1,0){60}}\put(160,40){\line(0,1){60}}
\put(100,55){\line(1,0){90}}\put(190,55){\line(0,1){45}}
\put(100,70){\line(1,0){105}}\put(205,70){\line(0,1){30}}
\put(100,85){\line(1,0){135}}\put(235,85){\line(0,1){15}}
\put(115,100){\line(0,-1){75}}\put(145,100){\line(0,-1){60}}
\put(175,100){\line(0,-1){45}}\put(220,100){\line(0,-1){15}}
\put(90,90){{\scriptsize9}}\put(90,75){{\scriptsize7}}
\put(90,60){{\scriptsize6}}\put(90,45){{\scriptsize4}}
\put(90,30){{\scriptsize2}} \put(150,-10){$\mu$}
\put(-100,-10){$\lambda$}
\end{picture}
\end{center}
\caption{Sylvester's bijection.}
\end{figure}
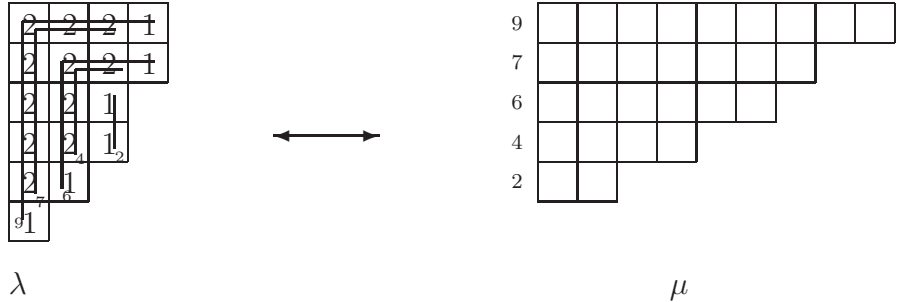

Applying the above bijection, we get the following correspondence between the set $\mathcal{T}(3)$ and the set $\mathcal{S}(3)$.
\[\begin{array}{lll}
3_r \leftrightarrows 3+3\ &
3_g \leftrightarrows 3\ & 2_r+1_r\leftrightarrows
2+2+1+1 \ \\[10pt]
 2_g+1_r\leftrightarrows
3+1\ & 2_r+1_g\leftrightarrows
3+2\ &
2_g+1_g\leftrightarrows
2+1+1\
\end{array}
\]
\[
\begin{array}{ll}  1_r+1_r+1_r\leftrightarrows
1+1+1+1+1+1\ &1_r+1_r+1_g\leftrightarrows
2+1+1+1\ \\[10pt]
1_r+1_g+1_g\leftrightarrows
2+2+1 \ &  1_g+1_g+1_g\leftrightarrows
1+1+1+1+1   \end{array}
\]

The following result immediately follows  from the combinatorial construction of Theorem \ref{main-AP}.

\begin{theo} Assume that $n\geq 1$, $r,l,p,q\geq 1$.
Let $s_{r,l,p,q}(n)$ denote the number of partitions $a_1+a_2+a_3+\cdots+a_{2\max\{r,l\}}$ satisfying $p+q\geq a_1\geq a_2
 \geq \cdots \geq a_{2\max\{r,l\}}\geq 0$ and $n=a_1+a_3+a_5+\cdots+a_{2\max\{r,l\}-1}$. Let $t_{r,l,p,q}(n)$ denote the number of two-color partitions of $n$ such that there are $r$ red parts and $l$ blue parts with the largest red part being not bigger than $p$ and the largest blue part being not bigger than $q$.
Then
\begin{equation*}
 s_{r,l,p,q}(n)=t_{r,l,p,q}(n)
 \end{equation*}
\end{theo}

\vspace{6mm}

\noindent{\bf Acknowledgments.} This work was supported by the National Science Foundation of China. We wish to thank the referees for valuable suggestions.


\begin{thebibliography}{99} \small

 \setlength{\itemsep}{-.8mm}

 \bibitem{Andrews-1976} G. E. Andrews, The Theory of Partitions,
 Addison-Wesley Publishing Co., 1976.

 \bibitem{and84} G.E. Andrews, Generalized Frobenius partitions,
Mem.  Amer. Math. Soc.  49 (1984), No. 301, iv+, 44 pp.

 \bibitem{Andrews-Paule-2021}  G. E. Andrews and P. Paule, MacMahon's partition analysis XIII: Schmidt type partitions and modular forms, J. Number Theory, (2021), DOI: 10.1016/j.jnt.2021.09.008.

 \bibitem{bes94}C. Bessenrodt, A bijection for Lebesgue's partition
identity in the spirit of Sylvester, Discrete Math. 132 (1994)
1--10.

\bibitem{bre99} D. Bressoud, Proofs and Confirmations, The story
of the alternating sign matrix conjecture, Cambridge University
Press, 1999.





\bibitem{mac84}P. A. Macmahon, Combinatory Analysis, Vol. II,
Cambridge University Press, Cambridge, 1915-1916, Reprinted:
Chelsea, New York, 1960.

\bibitem{yee03}A.J. Yee, Combinatorial proofs of generating
function identities for $F$-partitions, J. Combin. Theory Ser. A
102 (2003), 217--228.


 \bibitem{Wright-1965}E.M. Wright, An enumerative proof of an identity of Jacobi, J. London Math. Soc. 40 (1965)  55--57.



 \end{thebibliography}
\end{document}